\documentclass{amsproc}%
\usepackage{graphicx}
\usepackage{amscd}
\usepackage{amsmath}
\usepackage{amsfonts}
\usepackage{amssymb}%
\theoremstyle{plain}

\numberwithin{equation}{section}

\begin{document}
\title[Cayley Hamilton theorem with sandwich coefficients]{Cayley Hamilton theorem with sandwich coefficients for $n\times n$ matrices
over a ring satisfying $[x.y][u,v]=0$ }
\author{Jen\H{o} Szigeti}
\address{Institute of Mathematics, University of Miskolc, Miskolc, Hungary 3515}
\email{jeno.szigeti@uni-miskolc.hu}
\thanks{This research was carried out as part of the TAMOP-4.2.1.B-10/2/KONV-2010-0001
project with support by the European Union, co-financed by the European Social Fund.}
\subjclass[2010]{ 15A15,15A24,15A33,16S50}
\keywords{the symmetric determinant and characteristic polynomial, matrices with
commutator entries, Cayley-Hamilton theorem with matrix coefficients}

\begin{abstract}
If $A$ is an $n\times n$ matrix over a ring $R$ satisfying the polynomial
identity $[x,y][u,v]=0$, then an invariant Cayley-Hamilton identity of the
form%
\[
\underset{0\leq i,j\leq n}{\sum}A^{i}c_{i,j}A^{j}=0
\]
with $c_{i,j}\in R$ and $c_{n,n}=(n!)^{2}$ holds for $A$.

\end{abstract}
\maketitle

\noindent1. INTRODUCTION

\bigskip

The Cayley-Hamilton theorem and the corresponding trace identity play a
fundamental role in proving classical results about the polynomial and trace
identities of the $n\times n$ matrix algebra $\mathrm{M}_{n}(K)$ over a field
$K$ (see [2] and [3]).

In case of $\mathrm{char}(K)=0$, Kemer's pioneering work (see [5]) on the
T-ideals of associative algebras revealed the importance of the identities
satisfied by the $n\times n$ matrices over the Grassmann (exterior) algebra%
\[
E=K\left\langle v_{1},v_{2},...,v_{r},...\mid v_{i}v_{j}+v_{j}v_{i}=0\text{
for all }1\leq i\leq j\right\rangle
\]
generated by the infinite sequence of anticommutative indeterminates
$(v_{i})_{i\geq1}$. Accordingly, the importance of matrices over
non-commutative rings is an evidence in the theory of PI-rings, nevertheless
this fact has been obvious for a long time in other branches of algebra (e.g.
in the structure theory of semisimple rings). Thus a Cayley-Hamilton type
identity for such matrices seems to be of general interest.

In the general case (when $R$ is an arbitrary non-commutative ring with $1$)
Par\'{e} and Schelter proved (see [8]) that any matrix $A\in\mathrm{M}_{n}(R)$
satisfies a monic identity in which the leading term is $A^{k}$ for some
integer $k\geq2^{2^{n-1}}$\ and the other summands are of the form
$r_{0}Ar_{1}Ar_{2}\cdots r_{l-1}Ar_{l}$ with left scalar coefficient $r_{0}\in
R$, right scalar coefficient $r_{l}\in R$ and sandwich scalar coefficients
$r_{2},\ldots,r_{l-1}\in R$. An explicit monic identity for $2\times2$
matrices arising from the argument of [8] was given by Robson in [11]. Further
results in this direction can be found in [9] and [10].

Obviously, imposing extra algebraic conditions on the base ring $R$, we can
expect \textquotedblleft stronger\textquotedblright\ identities in
$\mathrm{M}_{n}(R)$. A number of examples show that certain polynomial
identities satisfied by $R$ can lead to \textquotedblleft
canonical\textquotedblright\ constructions providing an invariant
Cayley-Hamilton identity for $A$ of much lower degree than $2^{2^{n-1}}$.

If $R$ satisfies the polynomial identity%
\[
\lbrack\lbrack\lbrack...[[x_{1},x_{2}],x_{3}],...],x_{m}],x_{m+1}]=0
\]
of Lie-nilpotency (with $[x,y]=xy-yx$), then a left (and right)
Cayley-Hamilton identity of degree $n^{m}$ was constructed in [12]. Since $E$
is Lie-nilpotent of index $m=2$, this identity for a matrix $A\in
\mathrm{M}_{n}(E)$ is of degree $n^{2}$.

In [1] Domokos considered a slightly modified version of the mentioned
identity in which the left (as well as the right) coefficients are invariant
under the conjugate action of $\mathrm{GL}_{n}(K)$ on $\mathrm{M}_{n}(E)$. For
a $2\times2$\ matrix $A\in\mathrm{M}_{2}(E)$ the left scalar coefficients of
his Cayley-Hamilton identity are expressed as polynomials (over $K$) of the
traces $\mathrm{tr}(A)$, $\mathrm{tr}(A^{2})$ and $\mathrm{tr}(A^{3})$.

If $\frac{1}{2}\in R$ and $R$ satisfies the so called weak Lie-solvability%
\[
\lbrack\lbrack x,y],[x,z]]=0,
\]
then for a $2\times2$\ matrix $A\in\mathrm{M}_{2}(R)$ a Cayley-Hamilton trace
identity (of degree $4$ in $A$) with sandwich coefficients was exhibited in
[7]. If $R$ satisfies the identity%
\[
\lbrack x_{1},x_{2},...,x_{2^{t}}]_{\text{solv}}=0
\]
of general Lie-solvability, then a recursive construction (also in [7]) gives
a similar Cayley-Hamilton trace identity for $A\in\mathrm{M}_{2}(R)$ (its
degree depends on $t$).

In the present paper we consider an $n\times n$ matrix $A\in\mathrm{M}_{n}(R)$
over a ring $R$ (with $1$) satisfying the identity $[x,y][u,v]=0$ and
construct an invariant Cayley-Hamilton identity of the form%
\[
\underset{0\leq i,j\leq n}{\sum}A^{i}c_{i,j}A^{j}=0,
\]
where $c_{i,j}\in R$ are the sandwich coefficients and $c_{n,n}=(n!)^{2}$ is
the (central) leading coefficient.

We note that $[x,y][u,v]=0$ is the generating identity of the algebra
$\mathrm{U}_{2}(K)$\ of $2\times2$ upper triangular matrices (see [6]). The
identity $[x,y][x,z]=0$ (as well as $[[x,y],[x,z]]=0$) is a consequence of the
Lie-nilpotency $[[x,y],z]=0$ (see [4]). Clearly, the algebra $E$ shows that
$[x,y][u,v]=0$ is not a consequence of $[[x,y],z]=0$ and the algebra
$\mathrm{U}_{2}(K)$\ shows that $[[x,y],z]=0$ is not a consequence of
$[x,y][u,v]=~0$. Results about the logical relationships among the identities%
\[
\lbrack x,y][u,v]=0,[[x,y],z]=0\text{ and }[[x,y],[u,v]]=0
\]
can be found in [7].

We shall make extensive use of the so called symmetric characteristic
polynomial and the results in [12] and [13]. In order to provide a self
contained treatment in Section 2 we present all the necessary prerequisites.

\bigskip

\noindent2. CAYLEY-HAMILTON\ IDENTITY\ WITH\ MATRIX\ COEFFICIENTS

\bigskip

Let $R$ be an arbitrary ring with $1$. The preadjoint of a matrix
$A=[a_{i,j}]$ in $\mathrm{M}_{n}(R)$ was defined in [12] as $A^{\ast}%
=[a_{r,s}^{\ast}]$, where%
\[
a_{r,s}^{\ast}=\underset{\tau,\rho}{\sum}\mathrm{sgn}(\rho)a_{\tau
(1),\rho(\tau(1))}\cdots a_{\tau(s-1),\rho(\tau(s-1))}a_{\tau(s+1),\rho
(\tau(s+1))}\cdots a_{\tau(n),\rho(\tau(n))}%
\]
and the sum is taken over all permutations $\tau,\rho\in\mathrm{S}_{n}$ of the
set $\{1,2,\ldots,n\}$ with $\tau(s)=s$ and $\rho(s)=r$. The left and right
determinants of $A$ were defined in [13] as follows:%
\[
\mathrm{l\!\det}(A)=\mathrm{tr}(A^{\ast}A)\text{ and }\mathrm{r\!\det
}(A)=\mathrm{tr}(AA^{\ast}).
\]
If the base ring $R$ is commutative, then $\mathrm{tr}(AB)=\mathrm{tr}%
(BA)$\ for all $A,B\in\mathrm{M}_{n}(R)$. In spite of the fact that this well
known trace identity is no longer valid for matrices over a non-commutative
ring, the left and right determinants of $A$\ coincide (it was not recognized
in [13]).

\bigskip

\noindent\textbf{Proposition 2.1.}\textit{ The traces of the product matrices
}$A^{\ast}A$\textit{ and }$AA^{\ast}$\textit{ are equal: }$\mathrm{tr}%
(A^{\ast}A)=\mathrm{tr}(AA^{\ast})$\textit{.}

\bigskip

\noindent\textbf{Proof.} The trace of a matrix is the sum of the diagonal
entries, hence%
\[
\mathrm{tr}(A^{\ast}A)=\underset{1\leq r,s\leq n}{\sum}a_{r,s}^{\ast}%
a_{s,r}=\underset{\rho\in\mathrm{S}_{n},(\tau,s)\in\mathrm{S}_{n}^{\ast}}%
{\sum}\mathrm{sgn}(\rho)u(\rho,\tau,s),
\]
where $\mathrm{S}_{n}^{\ast}=\{(\tau,s)\mid\tau\in\mathrm{S}_{n}$, $1\leq
s\leq n$ and $\tau(s)=s\}$ and%
\[
u(\rho,\tau,s)=a_{\tau(1),\rho(\tau(1))}\cdots a_{\tau(s-1),\rho(\tau
(s-1))}a_{\tau(s+1),\rho(\tau(s+1))}\cdots a_{\tau(n),\rho(\tau(n))}%
a_{\tau(s),\rho(\tau(s))}.
\]
Similarly,%
\[
\mathrm{tr}(AA^{\ast})=\underset{1\leq r,p\leq n}{\sum}a_{p,r}a_{r,p}^{\ast
}=\underset{\rho\in\mathrm{S}_{n},(\alpha,p)\in\mathrm{S}_{n}^{\ast}}{\sum
}\mathrm{sgn}(\rho)v(\rho,\alpha,p),
\]
where%
\[
v(\rho,\alpha,p)\!=\!a_{\alpha(p),\rho(\alpha(p))}a_{\alpha(1),\rho
(\alpha(1))}\!\cdots\!a_{\alpha(p-1),\rho(\alpha(p-1))}a_{\alpha
(p+1),\rho(\alpha(p+1))}\!\cdots\!a_{\alpha(n),\rho(\alpha(n))}.
\]
Consider the following $\mathrm{S}_{n}^{\ast}\longrightarrow\mathrm{S}%
_{n}^{\ast}$ maps%
\[
(\tau,s)\longmapsto(\Theta(\tau,s),\tau(1))\text{ and }(\alpha,p)\longmapsto
(\Delta(\alpha,p),\alpha(n)),
\]
where the permutations $\Theta(\tau,s)$ and $\Delta(\alpha,p)$ in
$\mathrm{S}_{n}$ are defined by%
\[
(\!\tau(1)\!,\!1,2,\!\ldots\!,\!\tau(1)\!-\!1\!,\!\tau(1)\!+\!1\!,\!\ldots
\!,\!n\!-\!1,\!n)\!\overset{\Theta(\tau,s)}{\longmapsto}\!(\!\tau
(1)\!,\!\tau(2)\!,\!\ldots\!,\!\tau(s\!-\!1)\!,\!\tau(s\!+\!1)\!,\!\ldots
\!,\!\tau(n)\!,\!s)
\]
and%
\[
(1,2,\ldots,\alpha(n)\!-\!1,\alpha(n)\!+\!1,\!\ldots
\!,\!n\!-\!1\!,\!n\!,\!\alpha(n))\!\overset{\Delta(\alpha,p)}{\longmapsto
}\!(p\!,\!\alpha(1)\!,\!\ldots\!,\!\alpha(p\!-\!1)\!,\!\alpha
(p\!+\!1)\!,\!\ldots\!,\!\alpha(n)),
\]
respectively. It is straightforward to see that the above maps are mutual
inverses of each other:%
\[
\Delta(\Theta(\tau,s),\tau(1))=\tau,\Theta(\tau,s)(n)=s\text{ and }%
\Theta(\Delta(\alpha,p),\alpha(n))=\alpha,\Delta(\alpha,p)(1)=p.
\]
Since%
\[
u(\rho,\tau,s)=v(\rho,\Theta(\tau,s),\tau(1))\text{ and }v(\rho,\alpha
,p)=u(\rho,\Delta(\alpha,p),\alpha(n)),
\]
our claim is proved. $\square$

\bigskip

In view of Proposition 2.1,%
\[
\mathrm{s\!}\det(A)=\mathrm{tr}(A^{\ast}A)=\mathrm{tr}(AA^{\ast})
\]
can be called the \textit{symmetric determinant of }$A$. Let $R[x]$ denote the
ring of polynomials of the single commuting indeterminate $x$, with
coefficients in $R$. Let $[R,R]$ denote the additive subgroup of $R$ generated
by all commutators $[r,s]=rs-sr$ with $r,s\in R$. Using the unit matrix
$I\in\mathrm{M}_{n}(R)$, the \textit{symmetric characteristic polynomial} of
$A$ is the symmetric determinant of the $n\times n$ matrix $xI-A$ in
$\mathrm{M}_{n}(R[x])$:%
\[
p(x)=\mathrm{s\!}\det(xI-A)=\mathrm{tr}((xI-A)(xI-A)^{\ast})=\mathrm{tr}%
((xI-A)^{\ast}(xI-A)).
\]
The proof of Theorem 3.1 is based on the use of the following result from [13].

\bigskip

\noindent\textbf{Theorem 2.2.}\textit{ The symmetric characteristic polynomial
}$p(x)\in R[x]$\textit{ of a matrix }$A\in\mathrm{M}_{n}(R)$\textit{\ is of
the form}%
\[
p(x)=\lambda_{0}+\lambda_{1}x+\cdots+\lambda_{n-1}x^{n-1}+\lambda_{n}x^{n}%
\]
\textit{with }$\lambda_{0},\lambda_{1},\ldots,\lambda_{n-1},\lambda_{n}\in
R$\textit{ and }$\lambda_{n}=n!$\textit{. The product matrices }%
$n(xI-A)(xI-A)^{\ast}$\textit{ and }$n(xI-A)^{\ast}(xI-A)$\textit{ can be
written as}%
\[
n(xI-A)(xI-A)^{\ast}=p(x)I+C_{0}+C_{1}x+\cdots+C_{n}x^{n}%
\]
\textit{and}%
\[
n(xI-A)^{\ast}(xI-A)=p(x)I+D_{0}+D_{1}x+\cdots+D_{n}x^{n},
\]
\textit{where the matrices }$C_{i},D_{i}\in\mathrm{M}_{n}(R)$\textit{, }$0\leq
i\leq n$\textit{ are uniquely determined by }$A$\textit{. The entries of the
matrices }$C_{i},D_{i}$\textit{ are in }$[R,R]$\textit{, i.e. }$C_{i},D_{i}%
\in\mathrm{M}_{n}([R,R])$\textit{ for all }$0\leq i\leq n$. \textit{The right}%
\[
(\lambda_{0}I+C_{0})+A(\lambda_{1}I+C_{1})+\cdots+A^{n-1}(\lambda
_{n-1}I+C_{n-1})+A^{n}(n!I+C_{n})=0
\]
\textit{and the left}%
\[
(\lambda_{0}I+D_{0})+(\lambda_{1}I+D_{1})A+\cdots+(\lambda_{n-1}%
I+D_{n-1})A^{n-1}+(n!I+D_{n})A^{n}=0
\]
\textit{Cayley-Hamilton identities hold for }$A$\textit{.}

\bigskip

\noindent3. MATRICES OVER A\ RING\ SATISFYING $[x,y][u,v]=0$

\bigskip

\noindent\textbf{Theorem 3.1.}\textit{ If }$A\in\mathrm{M}_{n}(R)$\textit{ is
a matrix over a ring }$R$\textit{ satisfying the polynomial identity
}$[x,y][u,v]=0$\textit{, then an invariant Cayley-Hamilton identity of the
form}%
\[
\underset{0\leq i,j\leq n}{\sum}A^{i}c_{i,j}A^{j}=0
\]
\textit{holds for }$A$\textit{. The sandwich coefficients can be obtained as
}$c_{i,j}=\lambda_{i}\lambda_{j}$\textit{, where }$p(x)=\lambda_{0}%
+\lambda_{1}x+\cdots+\lambda_{n-1}x^{n-1}+\lambda_{n}x^{n}$\textit{ (with
}$\lambda_{n}=n!$\textit{) is the symmetric characteristic polynomial of }%
$A$\textit{.}

\bigskip

\noindent\textbf{Proof.} Rearranging the left and the right Cayley Hamilton
identities in Theorem 2.2, we obtain%
\[
\lambda_{0}I+A\lambda_{1}+\cdots+A^{n-1}\lambda_{n-1}+A^{n}\lambda_{n}%
=-(C_{0}+AC_{1}+\cdots+A^{n-1}C_{n-1}+A^{n}C_{n})
\]
and%
\[
\lambda_{0}I+\lambda_{1}A+\cdots+\lambda_{n-1}A^{n-1}+\lambda_{n}A^{n}%
=-(D_{0}+D_{1}A+\cdots+D_{n-1}A^{n-1}+D_{n}A^{n}).
\]
The multiplication of the above identities gives that%
\[
\underset{0\leq i,j\leq n}{\sum}A^{i}\lambda_{i}\lambda_{j}A^{j}%
=\underset{0\leq i,j\leq n}{\sum}A^{i}C_{i}D_{j}A^{j}.
\]
Now $C_{i}D_{j}=0$ is a consequence of $C_{i},D_{j}\in\mathrm{M}_{n}([R,R])$
and of $[x,y][u,v]=0$ in $R$. To complete the proof it is enough to note that
the coefficients $\lambda_{i}$, $0\leq i\leq n-1$ of the symmetric
characteristic polynomial of $A$\ are invariant under the conjugate action of
$\mathrm{GL}_{n}(\mathrm{Z}(R))$ on $\mathrm{M}_{n}(R)$, where $\mathrm{Z}(R)$
denotes the centre of $R$ (see [1] and [13]). $\square$

\bigskip

\noindent REFERENCES

\bigskip

\begin{enumerate}
\item M. Domokos, \textit{Cayley-Hamilton theorem for }$2\times2$%
\textit{\ matrices over the Grassmann algebra}, J. Pure Appl. Algebra 133
(1998), 69-81.

\item V. Drensky, \textit{Free Algebras and PI-Algebras}, Springer-Verlag, 2000.

\item V. Drensky and E. Formanek, \textit{Polynomial Identity Rings},
Birkh\"{a}user-Verlag, 2004.

\item S. A. Jennings, \textit{On rings whose associated Lie rings are
nilpotent}, Bull. Amer. Math. Soc. 53 (1947), 593-597.

\item A.R. Kemer, \textit{Ideals of Identities of Associative Algebras},
Translations of Math. Monographs, Vol. 87 (1991), AMS Providence, Rhode Island.

\item Yu.N. Malcev, \textit{A basis for the identities of the algebra of upper
triangular matrices}, (Russian), Algebra i Logika 10 (1971), 393-400; English
translation: Algebra and Logic 10 (1971).

\item J. Meyer, J. Szigeti and L. van Wyk, \textit{A Cayley-Hamilton trace
identity for }$2\times2$\textit{ matrices over Lie-solvable rings}, submitted

\item R. Par\'{e} and W. Schelter, \textit{Finite extensions are integral}, J.
Algebra 53 (1978), 477-479.

\item K.R. Pearson, \textit{A lower bound for the degree of polynomials
satisfied by matrices}, J. Aust. Math. Soc. Ser. A 27 (1979), 430-436.

\item K.R. Pearson, \textit{Degree 7 monic polynomials satisfied by a 3x3
matrix over a noncommutative ring}, Commun. Algebra 10 (1982), 2043-2073.

\item J.C. Robson, \textit{Polynomials satisfied by matrices}, J. Algebra 55
(1978), 509-520.

\item J. Szigeti, \textit{New determinants and the Cayley-Hamilton theorem for
matrices over Lie nilpotent rings}, Proc. Amer. Math. Soc. 125 (1997), 2245-2254.

\item J. Szigeti, \textit{Cayley-Hamilton theorem for matrices over an
arbitrary ring}, Serdica Math. J. 32 (2006), 269-276.
\end{enumerate}

\end{document}